\documentclass[10pt, a4paper, twocolumn]{article} 
\usepackage{graphicx}
\usepackage{subcaption}

\usepackage[english]{babel} 

\usepackage{microtype} 

\usepackage{amsmath,amsfonts,amsthm} 

\usepackage[svgnames]{xcolor} 

\usepackage{booktabs} 

\usepackage{lastpage} 

\usepackage{graphicx} 

\usepackage{enumitem} 
\setlist{noitemsep} 

\usepackage{sectsty} 
\allsectionsfont{\usefont{OT1}{phv}{b}{n}} 

\usepackage{geometry} 

\usepackage[T1]{fontenc} 
\usepackage[utf8]{inputenc} 

\usepackage{XCharter} 

\usepackage{fancyhdr} 
\fancypagestyle{firstpage}{ 
	\fancyhf{}
}

\newcommand{\institution}[1]{{\footnotesize\usefont{OT1}{phv}{m}{sl}\color{Black}#1}} 

\usepackage{titling} 

\newcommand{\HorRule}{\color{DarkGoldenrod}\rule{\linewidth}{1pt}} 

\pretitle{
	\vspace{-30pt} 
	\HorRule\vspace{10pt} 
	\fontsize{32}{36}\usefont{OT1}{phv}{b}{n}\selectfont 
	\color{DarkRed} 
}

\posttitle{\par\vskip 15pt} 

\preauthor{} 

\postauthor{ 
	\vspace{10pt} 
	\par\HorRule 
	\vspace{20pt} 
}

\usepackage{lettrine} 
\usepackage{fix-cm}	

\usepackage{xstring} 

\newcommand{\AuthorNames}{T. Ruggeri}
\newcommand{\ShortTitle}{In Memory of S. K. Godunov}

\title{In Memory of \\  Sergey K. Godunov 
\\	Personal Recollections} 

\author{ Tommaso Ruggeri 
	\newline\newline 
	 \institution{Dipartimento di Matematica, Università di Bologna, 40123 BOLOGNA, Italy, }\\
	   \institution{and Accademia Nazionale dei Lincei, ROMA, Italy,  tommaso.ruggeri@unibo.it}}

\date{} 

\usepackage{fancyhdr}
\begin{document}

\maketitle

\thispagestyle{firstpage}




\begin{center}
\begin{minipage}{0.95\columnwidth}
\noindent\textbf{\small{ABSTRACT :}}  
{\small{In this brief note, there is a short recollection of my scientific interactions with the great Russian mathematician Sergey Konstantinovich Godunov.}}
\end{minipage}
\end{center}

The first time I encountered the name of Sergey Konstantinovich Godunov was, of course, through his writings. Although I was not directly involved in numerical analysis at that time, it was universally known that Godunov was one of the greatest Russian mathematicians and mathematical physicists of the twentieth century. His seminal contributions to numerical methods for partial differential equations—especially the Godunov scheme for hyperbolic systems of conservation laws (1959)  \cite{God59}, had become fundamental pillars of numerical analysis and computational fluid dynamics, with countless applications. It was, and still is, a remarkable and extraordinarily influential achievement.

Around 1970 I began a collaboration with a young and brilliant Swiss–French researcher, Guy Boillat, who visited Bologna almost yearly to work with me. One day, Boillat came across a 1971 paper by the distinguished mathematicians Friedrichs and Lax, from the renowned Courant Institute in New York. That institute had become a refuge for many German mathematicians of Jewish origin who had emigrated to the United States during the Nazi period.

In their article, Friedrichs and Lax \cite{FL}  proved that every hyperbolic system of conservation laws compatible with a convex entropy could, after a suitable transformation, be rewritten in the so-called symmetric form—a notion introduced years earlier by Friedrichs himself. Boillat became intrigued by the question of whether one could find an explicit set of variables allowing the original system to be symmetrized.

In 1974, he indeed succeeded in finding such a change of variables \cite{Boillat}. The result immediately struck me as both elegant and important, and I asked him to give a seminar in Bologna. In the audience there was my colleague Mauro Fabrizio, who pointed out that Godunov had already obtained such a transformation in 1961 in a brief note, showing not only that the Euler system but all variational systems could be symmetrized \cite{Gods}. Although Boillat’s result was more general, he had to acknowledge that the idea itself was not entirely new. From that moment on, the terminology ``Godunov symmetric systems'' became standard in the literature.

I believe Godunov appreciated this recognition: in later years, especially in his works with his longtime collaborator Evgeny Romenski, Boillat is often cited (see, e.g. \cite{RG}).  I myself became interested in extending the method to the relativistic case, where Boillat’s approach could not be applied. Together with my first student, Alberto Strumia, we discovered that the appropriate symmetrizing variables were precisely the multipliers of the system yielding the supplementary balance law \cite{RS}. This established a natural bridge between my work and the Siberian school founded by Godunov.

Some years later—perhaps at the end of 1990s—I received an email from Sergey Gavrilyuk, who had recently moved to Marseille with the support of Henri Gouin, whom I did not yet know personally. Gavrilyuk proposed a European project involving Novosibirsk (represented by Godunov), Marseille, and Bologna. We drafted what I believed to be an extremely interesting and timely project, but it was rejected with the rather discouraging comment that it was ``too ambitious and too general.'' Nevertheless, we started collaborating anyway. Gouin invited me to Marseille as a visiting professor, and this was the beginning of a long and fruitful friendship with him and with Gavrilyuk; for many years I was a regular visitor to Marseille, and Gouin frequently visited Bologna as well.

In 1999, Sergey Godunov himself was invited to Marseille as a visiting professor. That was when I finally had the pleasure of meeting this extraordinary scientist in person. I found him a very kind man, not very tall, dressed in a sort of ``colonial-style'' jacket that, as I later learned, was typical of his personal taste. He delivered several seminars, and I was deeply fascinated by both his clarity and his intellectual energy.

There was, however, a small linguistic challenge: he spoke French—though he understood some English—while I understood French well but did not speak it fluently. As a result, our conversations often proceeded in a curious manner: I spoke to him in English and he answered in French, each of us understanding the other but replying in a different language. We often took long walks around the campus or went to lunch together, and gradually we became more familiar with each other.

One day, out of curiosity, I asked whether he had ever read the paper by Friedrichs and Lax. I should never have asked! He suddenly turned red and began speaking rather vehemently—I am still not quite sure in which language—explaining that he had been invited to the Courant Institute years before and had presented his results on symmetrization. Nevertheless, the authors later wrote their paper without citing him, and despite the methodological differences, they had completely ignored his prior work. He remained genuinely upset about this episode, and understandably so.

A more light-hearted recollection concerns our lunches with Henri and Sergey at a small traditional restaurant near the University, where the food was excellent. At the end of the meal, they always offered cheese, and Godunov absolutely loved it. While we would take a small slice, he would happily finish almost the entire tray!

Because of our growing collaborations, I invited Godunov—and of course Gavrilyuk and Gouin—to the international conference WASCOM (Waves and Stability in Continuous Media), which I organized every two years. He attended the 2001 edition in Porto Ercole and delivered a lecture in French. I remember vividly his scientific liveliness and also his somewhat polemical temperament: during a talk by Giovanni Russo, an excellent numerical analyst, he criticized the presentation rather strongly, arguing that he saw nothing new in it. Most likely it was a misunderstanding, perhaps due to language.

\begin{figure}[ht]
	\centering
	
	\begin{subfigure}{0.88\columnwidth}
		\centering
		\includegraphics[width=\columnwidth]{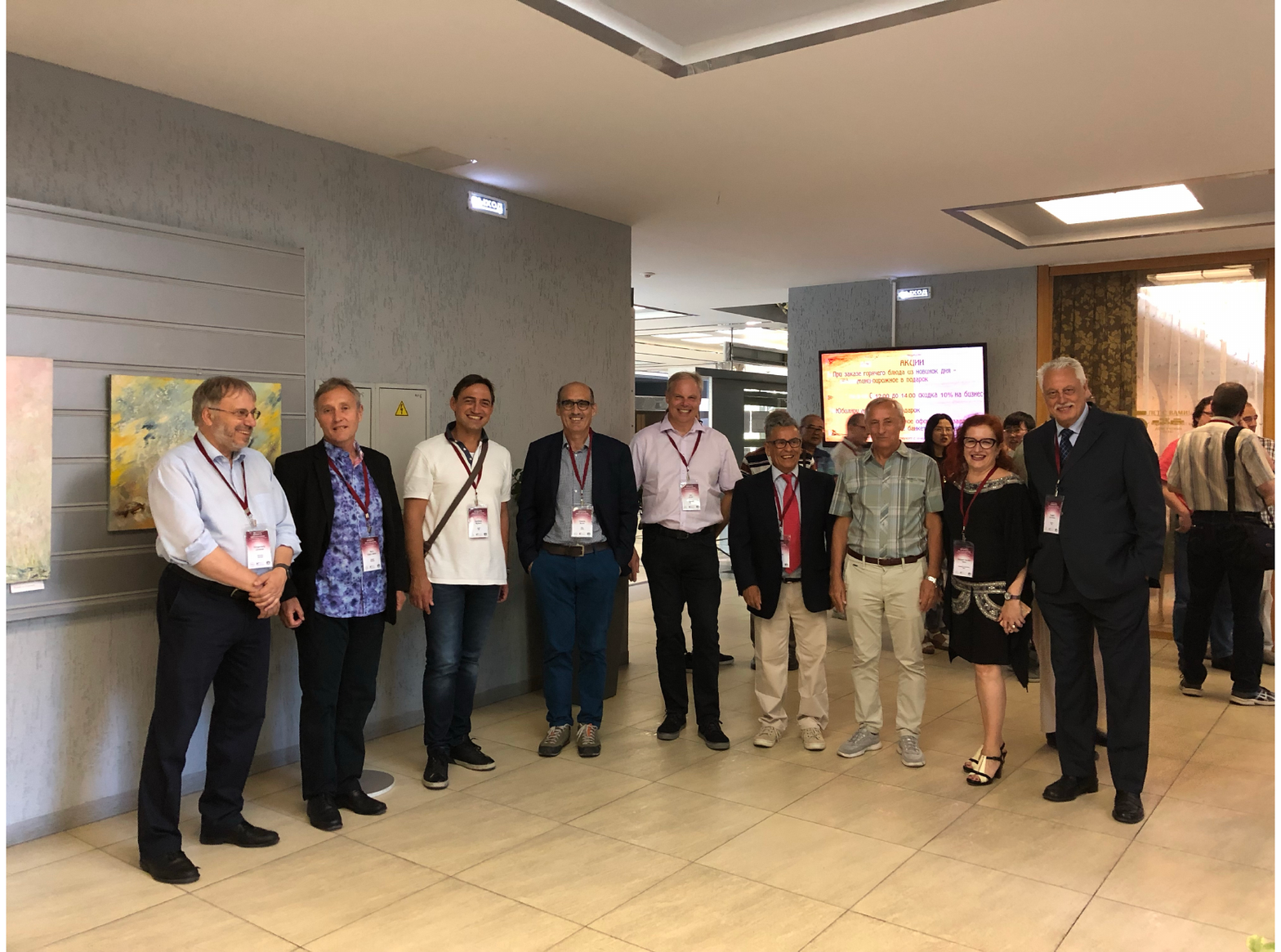}
		\caption{Group of participants. From left: Christian Klingenberg, 
			Claus Dieter Munz, Michael Dumbser, Bruno Després, Armin Iske, 
			Eleuterio Toro, Evgeny Romenski, Elena Vázquez Cendón and 
			Tommaso Ruggeri.}
	\end{subfigure}
	
	\hfill
	
	\begin{subfigure}{0.68\columnwidth}
		\centering
		\includegraphics[width=\columnwidth]{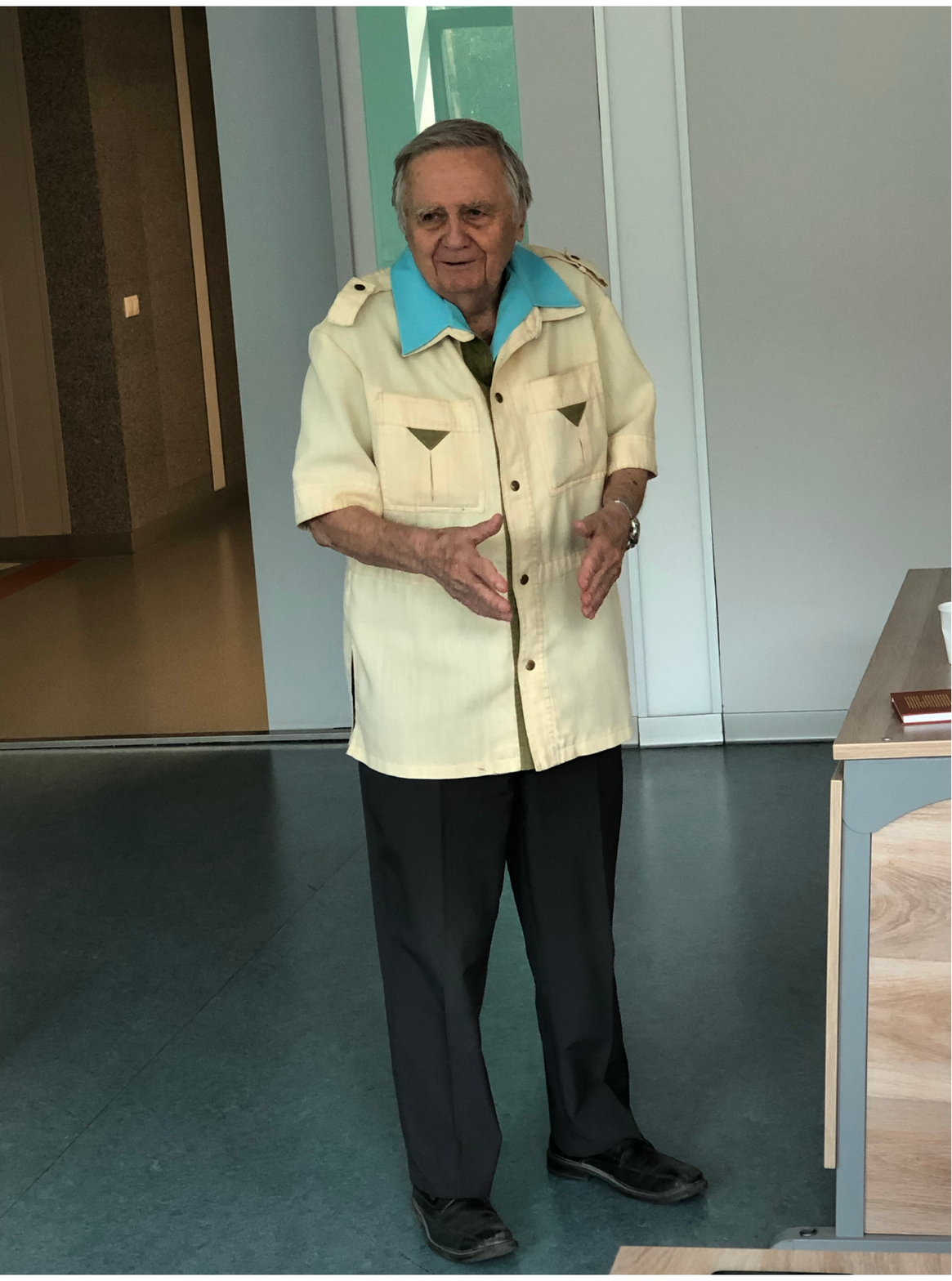}
		\caption{Prof.\ S.\ K.\ Godunov during the meeting celebrating 
			his 90th birthday.}
	\end{subfigure}

\end{figure}

Later, thanks to my friendship with Evgeny Romenski, for whom I have great admiration and affection, I had the honor of being invited to Novosibirsk, first for Godunov’s 80th birthday (an event I unfortunately could not attend) and later, in 2019, for his 90th birthday conference. That celebration was truly impressive. The opening session was held in a huge auditorium, and I immediately sensed how deeply respected and admired Godunov was in his own country. Politicians attended, and newspapers dedicated extensive coverage to him—recognition he fully deserved. Unlike most scientific conferences today, the official language of the opening day was Russian, with simultaneous translation.

I had the honor of giving one of the plenary lectures on that first day. There, I truly understood the extraordinary mathematical level of the Russian school, particularly the group at Novosibirsk and the celebrated Sobolev Institute, which in the late 1950s and early 1960s had gathered some of the finest mathematicians of the Soviet Union. As is well known, the creation of this scientific center—far from Moscow and Leningrad—was driven by the vision of Keldysh, Sobolev, Lavrent'ev and other eminent scientists, seeking not only strategic decentralization but also greater scientific freedom. Among the young talents they recruited to Siberia, Godunov was certainly one of the brightest. He found an ideal environment there and became a central figure of the Siberian school in numerical methods and hyperbolic equations.

I spent wonderful days in Novosibirsk, meeting old and new friends. I am not sure whether, in the bustle of the event, Godunov recognized me, and I was sorry that he felt tired and did not attend the gala dinner. That was the last time I saw him. 

Some Colleagues later put together a commemorative volume, and I contributed a paper dedicated to him \cite{book}. I enclose two  photographs were taken during the conference celebrating 
	the 90th birthday of Prof.\ S.\ K.\ Godunov, held in Novosibirsk in 
	August~2019.

I was deeply saddened when I learned of his passing. From what his students told me, he could indeed have a difficult character, but in the few occasions when I met him, I always found him pleasant, sincere, and unmistakably different from other scientists. I truly believe that his name will remain forever in the history of mathematics, and I am proud to have had the honor of knowing him personally and of having worked on topics in which he was a true pioneer.

\end{document}